\documentclass[twocolumn,showpacs,preprintnumbers,amsmath,amssymb,footinbib]{revtex4}

\usepackage{graphicx}
\usepackage{dcolumn}
\usepackage{bm}
\usepackage{url}

\begin{document}

\preprint{arXiv preprint}

\title{A method for dense packing discovery}

\author{Yoav Kallus}
 \email{yk328@cornell.edu}
\author{Veit Elser}%
\affiliation{%
Laboratory of Atomic and Solid-State Physics, Cornell University, Ithaca, New York, 14853
}%

\author{Simon Gravel}
\affiliation{
Department of Genetics, Stanford University School of Medicine, Stanford, California 94305-5120
}%

\date{\today}

\begin{abstract}
The problem of packing a system of particles as densely as possible is foundational in the field of
discrete geometry and is a powerful model in the material and biological sciences. As packing
problems retreat from the reach of solution by analytic constructions, the importance of an efficient
numerical method for conducting \textit{de novo} (from-scratch) searches for dense packings becomes crucial. In this paper, we use the \textit{divide and
concur} framework to develop a general search method for the solution of periodic constraint
problems, and we apply it to the discovery of dense periodic packings. An
important feature of the method is the integration of the unit cell parameters with the other
packing variables in the definition of the configuration space.
The method we present led to improvements in the
densest-known tetrahedron packing which are reported in Ref. \cite{tetra}.
Here, we use the method to reproduce the densest known lattice sphere packings and 
the best known lattice kissing arrangements in up to 14 and 11 dimensions respectively (the first
such numerical evidence for their optimality in some of these dimensions).
For non-spherical particles, we report a new
dense packing of regular four-dimensional simplices with density $\phi=128/219\approx0.5845$
and with a similar structure to the densest known tetrahedron packing.

\end{abstract}

\pacs{61.50.Ah, 45.70.-n, 02.70.-c}
\maketitle

\section{\label{sec1}Introduction}

The dense packing behavior of a general solid body (particle) in a Euclidean
space is a problem of interest in mathematics, physics,
and many other fields. A packing is a collection of
particles in the Euclidean space $\mathbb{R}^d$, wherein no two particles overlap
(i.e., the intersection of any two particles has an empty interior) and the packing fraction
or density $\phi$
is then the volume fraction of space covered by the particles.
Of particular interest are packings of a given particle (wherein all particles are congruent), 
and the problem of interest is to determine the maximum possible density
$\phi_\text{max}$ among all packings of a given particle. A packing that realizes this
maximum can be thought of as the equilibrium state of the system of classical hard
particles in the limit of infinite pressure or zero temperature.

The general problem of packing congruent particles was posed as a
part of the eighteenth of David Hilbert's famous \textit{Mathematische Probleme}:
\begin{quote} How can one arrange most densely in space an infinite number of
equal solids of a given form, e.g., spheres with given radii or regular
tetrahedra with given edges (or in prescribed position), that is, how can one so
fit them together that the ratio of the filled to the unfilled space may be as
large as possible? \cite{Hilbert} \end{quote}
This part of the problem has been
taken over the years  as the resolution of the Kepler conjecture
about the densest packing of spheres in three dimensions \cite{HilbertChallenge}, and has therefore
been considered resolved since the latter was proved by Hales \cite{Hales}. However,
Hilbert's statement
of the problem does not single out the sphere, and actually
mentions the regular tetrahedron as another particle of interest. Recent work
diverging from the focus on spherical particles has spotlighted ellipsoids
\cite{ellipse}, regular and semi-regular polyhedra \cite{TorqNat,TorqPRE} (and
the regular tetrahedron in particular \cite{Chen,Glotzer,tetra,TorqUniform,ChenUniform}), and
superballs \cite{superball}. Few bounds are known for the maximum packing
fraction of general convex particles. Kuperberg and Kuperberg have shown that
for any convex particle in two dimensions,
$\phi_\text{max}\ge\sqrt{3}/2\approx0.86602$ \cite{Kup1}. Torquato \textit{et al.} used
the known maximal packing density of spheres to derive an upper bound on the packing
density of any solid, but this bound is trivial (i.e., $\phi_\text{max}\le\phi^U$, where $\phi^U>1$)
for many solids \cite{TorqNat}. Ulam has conjectured
that in three dimensions, the sphere achieves the lowest maximum packing
fraction, $\phi_\text{max}=\pi/\sqrt{18}\approx0.74048$, among all convex
particles \cite{Ulam}.

In the quest for dense packings of various particles, analytic and numerical investigations have
both played important roles. The former have been very successful in the study of the dense packing of spheres, where analytic constructions based on groups, codes, and laminated
lattices have produced the densest-known sphere packings and lattice sphere packings in many
dimensions \cite{SPLAG}. However, the analytic approach to the construction of dense packings
relies on the imagination of the constructor, and for a variety of other problems the
densest packings have evaded the creativity of analytic investigators and were only uncovered
in computational investigations. While complete (i.e., exhaustive) algorithms exist for some problems (such as the algorithm
in Ref. \cite{Covering}, which gave new best known results for the lattice covering and covering-packing
problems in some dimensions), they do not exist or have runtimes that are too long for other
problems. In those cases, incomplete search algorithms become necessary.

One example of a dense packing that has only been uncovered by a \textit{de novo} 
numerical search
is the currently densest-known packing of tetrahedra, whose structure 
was first hinted at by a numerical search using the method described in this paper \cite{tetra}.
The structure was later optimized by Torquato and Jiao \cite{TorqUniform}
and by Chen \textit{et al.} \cite{ChenUniform}.
Results of subsequent Monte Carlo simulations have reproduced this structure and 
suggest it is the densest packing of regular
tetrahedra at least with a small number ($\le16$) of tetrahedra in the unit
cell \cite{TorqUniform,ChenUniform}.
Another \textit{de novo} search with Monte Carlo dynamics
has uncovered a packing based on a quasicrystal approximant reminiscent of
the Frank-Kasper $\sigma$-phase with a slightly lower density \cite{Glotzer}.
As these two structures were overlooked by previous
analytical investigations \cite{Chen, Welsh}, it is quite likely that without 
the results of \textit{de novo} searches, they would have remained unimagined and undiscovered.

In the best case, such searches would produce the optimal packing possible subject to the
built-in restrictions (such as number of particles in the unit cell or unit cell shape). However,
in problems exhibiting a large degree of frustration, the presence of many local optima
that are separated from each other by high barriers complicates the task of finding the optimal
packing. The tendency of simulations to get stuck in the local optima of  such a rugged optimization
landscape, especially when these local optima proliferate as more 
particles are simulated, has been held responsible for suboptimal results in searches \cite{TorqNat,
TorqPRE}. One technique which has been observed to relieve dynamical stagnation 
in Monte Carlo simulations at high pressures has been to allow slightly unphysical moves,
such as allowing particles to temporarily overlap \cite{Glotzer}. 

We propose a novel search method as an alternative to Monte Carlo simulations, with
a number of features that directly address these observations. The method 
is based on the dynamics of the difference map, a constraint-satisfaction
iterative search algorithm, and on the \textit{divide and concur} constraint framework
(we abbreviate this combination $D-C$, where the minus sign stands for the difference map)
\cite{DM-PNAS,D-C,Gravel}. It adapts the $D-C$ approach to the case of periodic problems
and we shall call it \textit{periodic divide and concur} (PDC).  
The difference map
is designed to avoid being trapped in local optima
and has been demonstrated in multiple applications to find solutions of
highly non-convex problems, including finite packing problems with large
numbers of particles, from random starting
configurations \cite{DM-PNAS,D-C,Gravel}. The search proceeds through a non-physical configuration
space, cutting through the conventional physical optimization landscape.
Still, it is to be expected that the exponential growth in the number 
of local optima in the configuration space, which the search will still have to traverse,
will nevertheless lead to suboptimal results when many independent particles are
included in the search. 
Therefore, as discussed below, it is crucial for the unit cell variables to be
aggressively optimized so that the number of particles to be simulated can be
reduced. The incorporation of the unit cell variables directly into the basic dynamics of the search
achieves this goal.

Numerical searches are restricted to finite-dimensional configuration spaces, and therefore have
been largely limited to investigating periodic packings, packings which are preserved under translations
by a lattice $\Lambda$. In a general periodic packing, the particles 
are partitioned into $p$ orbits of the lattice $\Lambda$, and when $p=1$ the packing
is called a lattice packing. In physics, any periodic arrangement is usually referred to as a lattice and
the special case of $p=1$ is known as a Bravais lattice.
In general, the maximum density need not be realizable by
a periodic packing, but arbitrarily close densities are realizable with
periodic packings of arbitrarily large $p$. Similarly, arbitrarily accurate approximations
of any packing can be obtained using a sufficiently large cubic or orthorhombic unit cell.
However, due to the rapid increase in computational complexity
and the proliferation of local optima as the number of independent particles rises,
it is often preferable to include fewer particles but allow for a variable unit cell shape.
We focus then on searching for packings with a small number of particles in the unit cell.

To our knowledge, variable unit cells have only been introduced recently to searches
for dense packings, for instance with the adaptive shrinking cell scheme in Refs.
\cite{ellipse,TorqPRE} and with the
use of Parrinello-Rahman dynamics in the space of lattices in Ref. \cite{Cohn}.
The increased particle population associated with restricting unit cell variability
can sometimes be tolerated in two and three dimensions,
but in high dimensions the number of particles that must be simulated grows exponentially
due to the curse of dimensionality and this approach becomes impractical. The
constraint-satisfaction formulation of the periodic packing problem
used in PDC features a variable unit cell and naturally treats
the positions of particles in the unit cell and the unit cell parameters on the same footing.
This new approach allows us to successfully look for dense
sphere packings in dimensions as high as 14, further than probed by any previously
reported unbiased numerical exploration of periodic packings.

Besides the density of a packing, another attribute of interest is the coordination number, that is,
the number of nearest neighbors of particles in the packing.
In the case of spherical particles, this amounts to the number of spheres in contact with a given sphere,
known as the kissing number \cite{SPLAG}.
Searching for high-coordination number arrangements around a single
sphere has been accomplished previously with the $D-C$ method \cite{D-C}.
Here we apply PDC to search
for space-filling periodic arrangements of high coordination number, and particularly lattice
arrangements.

An efficient \textit{de novo} numerical search method can provide critical utility in the field of packing. In addition to the ability of a \textit{de novo} search to
provide confidence in a putative, but not proven, optimal result, a \textit{de novo} search has often been responsible for surprising new results: two recent examples in which unexpected
(as it turns out, quasiperiodic) packings were found as the
results of \textit{de novo} searches are in the problem of
tetrahedron packing \cite{Glotzer} and in the ten-dimensional kissing number problem \cite{consQ}.
It is with these motivations that we introduce the PDC method in this paper. 
In Section \ref{sec1a} we introduce the $D-C$ scheme by presenting a simple example which
serves to motivate the constructions in the subsequent sections.
In Section \ref{sec2} we formulate the problems tackled in this paper
--- sphere packing, the lattice kissing number, and polytope packing ---
in terms of constraint satisfaction.
In Section \ref{sec3} we describe in detail aspects of our implementation of the PDC 
search, including efficient computation of projections to the
constraints of Section \ref{sec2}. In Section \ref{sec4} we present
some results of PDC for the problems discussed,
including a newly discovered packing of regular four-dimensional simplices.
In Section \ref{sec5} we present concluding remarks.

\section{\label{sec1a}Motivation}
\subsection{\label{sec1a1}The $D-C$ scheme}

The key step in applying the $D-C$ approach to packing problems is to recast the problem
as a problem of constraint satisfaction. Particularly, we must express it as the problem of
finding a configuration in a Euclidean configuration space ($\Omega$), which satisfies two constraints.
We identify a constraint $C$ with the subset $C\subseteq\Omega$ of configurations satisfying the
constraint. A projection of a configuration $x$ to a constraint $C$ is the operation of finding a configuration $x' \in C$ that minimizes the distance $||x-x'||$.
Each of the two constraints ($C, D \subseteq \Omega$), must be simple enough
that the operation of projecting an arbitrary configuration to it can be computed efficiently.
The iterative map used in exploring the configuration space takes advantage of the formulation
of the problem in terms of two simple constraints, as outlined in section \ref{sec1a2}.
In this section we present the application of the $D-C$ scheme to finite sphere packing
problems, which has been developed and implemented in Ref. \cite{D-C}, as an introduction
to the main ideas of the scheme.

\begin{figure}

\includegraphics[scale=0.36]{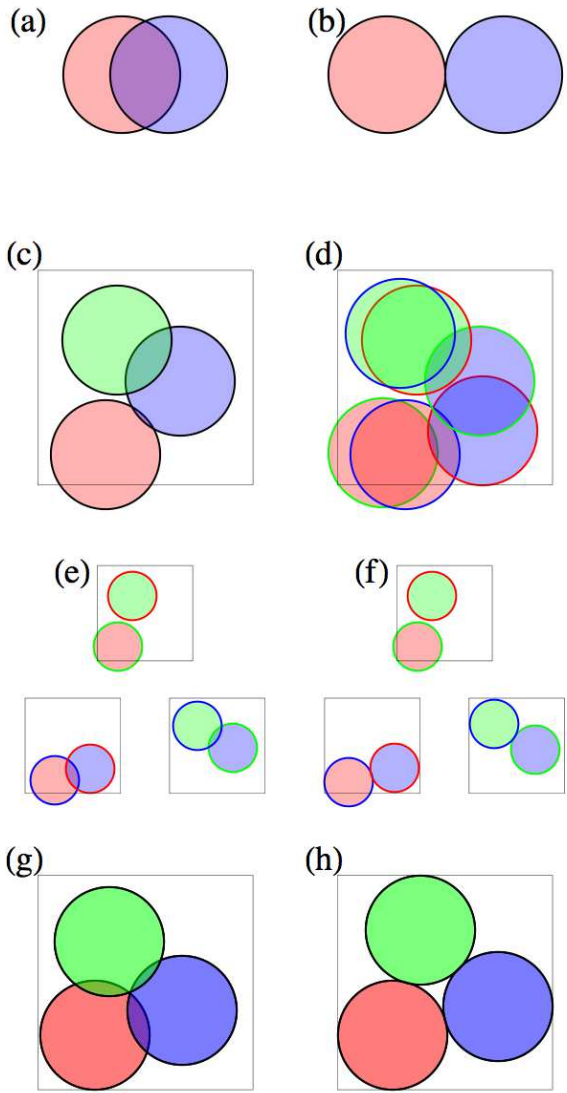}

\caption{
\label{fig0}
An illustration of the $D-C$ scheme in the case of packing three disks into a square box.
In the case of two overlapping disks (a), the exclusion constraint is simple in the sense
that a projection to the constraint can be performed efficiently. The projection is given by
\eqref{sphereproj} and yields the configuration (b). In the case of three disks (c),
there is no similarly efficient projection method to the constraint that no overlaps occur. In the $D-C$
scheme, each disk is represented by two replicas (d), which together make three independent
replica pairs (e). The exclusion constraint, now also called the ``divide'' constraint,
is modified so that only replica pairs are prohibited
from overlapping, and any other overlaps are allowed. Thus, the projection to the exclusion
constraint can be performed independently on each replica pair as in the case of two disks (f).
A second constraint, the ``concur'' constraint, requires all replicas representing
a single disk to coincide and requires the disk to lie within the confinement box.
The result of projecting the configuration (d) to this constraint is the configuration (g).
In order to search for a configuration satisfying both constraints, we do not alternately
project from one constraint to the other, but instead use the difference map \eqref{DMeq}
to evolve the configuration. The result of a successful search is a configuration (h) satisfying both
constraints, which by construction corresponds to a solution of the problem.}
 \end{figure}

The defining constraint of packing problems is the constraint that no
particles in the packing overlap, which we call the {\it exclusion constraint}. As a simple illustration
of this constraint, consider the exclusion of a pair of unit-radius disks in $\mathbb{R}^2$. In this case,
the configuration space $\Psi$ is parameterized by the positions of the centers of the two disks:
\begin{equation}
\Psi = \{ (\mathbf{x}_1,\mathbf{x}_2) \colon \mathbf{x}_1,\mathbf{x}_2 \in \mathbb{R}^2 \}\text.
\end{equation}
The exclusion constraint is then
\begin{equation}
K_\text{excl} = \{ (\mathbf{x}_1,\mathbf{x}_2) \in \Psi \colon ||\mathbf{x}_1 - \mathbf{x}_2|| \ge 2\} \subseteq \Psi\text.
\end{equation}
This constraint adheres to the simplicity criterion of having an efficient method to compute a
projection to it. Specifically, the projection is given by
\begin{equation}\label{sphereproj}
\pi_{K_\text{excl}} \left[ (\mathbf{x}_1,\mathbf{x}_2) \right] = \begin{cases} 
(\mathbf{x}'_1,\mathbf{x}'_2) & \text{if } ||\mathbf{x}_1 - \mathbf{x}_2|| < 2 \\
(\mathbf{x}_1,\mathbf{x}_2) & \text{otherwise,}
\end{cases}
\end{equation}
where,
\begin{align}
\mathbf{x}'_1 = & \mathbf{x}_1 + \frac{2-||\mathbf{x}_1 - \mathbf{x}_2||}{2||\mathbf{x}_1 - \mathbf{x}_2||}(\mathbf{x}_1 - \mathbf{x}_2) \\
\mathbf{x}'_2 = & \mathbf{x}_2 - \frac{2-||\mathbf{x}_1 - \mathbf{x}_2||}{2||\mathbf{x}_1 - \mathbf{x}_2||}(\mathbf{x}_1 - \mathbf{x}_2) 
\end{align}
as illustrated in Figure \ref{fig0}.

A more complicated case arises when three or more disks are considered. In this case, the exclusion
constraint,
\begin{align*}
K_\text{excl} &=  \{ (\mathbf{x}_1,\ldots \mathbf{x}_n) \in \Psi \colon\\
&||\mathbf{x}_i - \mathbf{x}_j|| \ge 2 \text{ for all } 1\le i < j \le n\}\text,
\end{align*}
is not a simple constraint according to the criterion above.
Alternatively, we could replace $K_\text{excl}$ by many pairwise exclusion constraints
\begin{equation}
K_\text{excl}^{i,j} = \{ (\mathbf{x}_1,\ldots \mathbf{x}_n) \in \Psi \colon ||\mathbf{x}_i - \mathbf{x}_j|| \ge 2\}\text.
\end{equation}
The pairwise constraints
are all individually simple. However, as noted above, we are limited to problems described
by only two simple constraints.

{\it Divide and concur} provides a general procedure for reducing
the number of simple constraints to two, at the expense of enlarging the configuration space.
This reduction is achieved
by parameterizing the configuration space with more variables than are necessary to fully
specify a configuration.
In the example at hand, the new configuration space is
\begin{equation}
\Omega = \{ (\mathbf{x}_{1,2},\ldots \mathbf{x}_{n,n-1}) \colon \mathbf{x}_{i,j} \in \mathbb{R}^2
\text{ for all } i\neq j\}\text,
\end{equation}
where we call all the variables $\mathbf{x}_{i,j}$ for a particular index $i$ the {\it replicas} of the original
variable $\mathbf{x}_{i}$. Every configuration $(\mathbf{x}_1,\ldots \mathbf{x}_n) \in \Psi$
can be identified with a configuration $(\mathbf{x}_{1,2},\ldots \mathbf{x}_{n,n-1}) \in \Omega$,
wherein $\mathbf{x}_{i,j}=\mathbf{x}_{i}$ for all $i,j$, through a simple linear map $A$.
Enough redundant variables have been introduced to the
configuration space so that each pairwise exclusion constraint can now be written 
in terms of a private set of variables, disjoint from the private variables of other constraints:
\begin{equation}
D^{i,j} = \{ (\mathbf{x}_{1,2},\ldots \mathbf{x}_{n,n-1}) \in \Omega \colon ||\mathbf{x}_{i,j} - \mathbf{x}_{j,i}|| \ge 2\}\text.
\end{equation}
The intersection, $D\subseteq \Omega$, of all of the pairwise exclusion constraints, which we will call
the ``divide'' constraint, is now also simple,
since the projection can be performed independently on each set of private variables (Figure \ref{fig0}).

The map $A:\Psi\to\Omega$ from the original configuration space (the {\it physical
configuration space}) to the new one (the {\it formal configuration space}) is not surjective,
and so a general point in the formal configuration space does
not correspond to a valid physical configuration. The ``concur'' constraint $C=A(\Psi)$ is given by the
range of $A$, the subset of $\Omega$ that does correspond to valid configurations.
That is, the constraint requires redundant specifications of an original variable to concur in regard
to its value. Since $A$ is linear, $C$ is also a simple constraint.

Another constraint that must usually be addressed in packing problems with a finite number
of particles is the confinement constraint. In most cases the particles, or their centers, are confined
to lie in some subset $M$ of space, where $M$ can be either some region of finite volume, or
a compact manifold (as in the case of spherical codes). As a subset of the original
configuration space, the confinement constraint is written as
\begin{equation}
K_\text{conf} = \{ (\mathbf{x}_{1},\ldots \mathbf{x}_{n}) \in \Psi \colon \mathbf{x}_{i} \in M \text{ for all }i\} \subseteq \Psi\text.
\end{equation}
We can incorporate this constraint into the ``concur'' constraint, $C$, by modifying it to be the image
$A(K_\text{conf})$ instead of the entire range of $A$. In our example, this would give the constraint
\begin{equation}
C = \{ (\mathbf{x}_{1,2},\ldots \mathbf{x}_{n,n-1}) \in \Omega \colon \mathbf{x}_{i,j} = \mathbf{x}_{i} \in M \text{ for all }i,j\}\text.
\end{equation}
Since $A$ is linear, the projection to $C=A(K_\text{conf})$ can be
decomposed into a projection to $A(\Psi)$ followed by a projection to $A(K_\text{conf})$.
The first step is performed by taking the average position of all the
replicas of each disk. The second step is performed by projecting this average position to $M$ (see
Figure \ref{fig0}). In general, this two-step projection method is valid for handling the constraints in the
physical configuration space that are simple at the outset and do not require the introduction of
new variables.

The result of the above construction is that a configuration in $\Omega$ satisfies the
``divide'' and ``concur'' constraints simultaneously if and only if it corresponds to a configuration
in $\Psi$ which satisfies all the exclusion constraints and the confinement constraint;
that is, it corresponds to a solution of the packing problem under consideration.

In the following sections we modify the above simple construction
so as to generalize the method in two major ways.
The first generalization is to packings of infinite regions, instead of only finite ones.
Specifically, we allow for periodic packings with an arbitrary unit cell. This is achieved by
generalizing the idea of replicas of particles to include also their periodic images.
When the unit cell vectors are included in the original set of parameters, the map
from the original parameter space
to the space of replica configurations is still linear, though a little more elaborate.
Additionally, the confinement constraint of finite packings is replaced in the case
of periodic packings by a constraint on the unit cell volume, ensuring a specified density.

The second generalization is to packings of non-spherical particles, specifically convex
polytopes. This is achieved by representing each particle not only by the position of its
centroid, but by the positions of all its vertices. A new constraint, the rigidity constraint,
is added to ensure that the particle is not deformed in the solution. Despite the mathematical
complications that arise from these two generalizations, the conceptual framework is
identical to the above example, and the constructions in the
following sections will draw attention to the analogy with the construction presented above. 

\subsection{\label{sec1a2}The difference map}

Given a problem formulated as the task of finding a configuration $\mathbf{x}\in C\cap D$,
simultaneously satisfying the constraints $C,D\subseteq \Omega$, we wish to use
the availability of efficient methods for computing the projections $\pi_C$ and $\pi_D$ to
the constraints in order to set up an iterated map to search through the configuration
space for a solution. Naive schemes, such as the alternating projections map $\mathbf{x}\mapsto
\pi_D ( \pi_C(\mathbf{x}) )$, suffer from the problem of stagnation at near solutions (local minima
of the distance between the two constraints). The difference map, a slightly more sophisticated
scheme, is designed to provide efficient search dynamics while avoiding the traps of
local minima \cite{DM-PNAS}.

The difference map (DM) can be written in terms of the projections $\pi_C$ and $\pi_D$
and one parameter $\beta$:
\begin{equation}
\operatorname{DM}: \Omega \to \Omega
\end{equation}
\begin{equation}\label{DMeq}
\mathbf{x} \mapsto \mathbf{x} + \beta 
\left[ \pi_D\left( f_C(\mathbf{x})\right)
-\pi_C\left( f_D(\mathbf{x})\right)\right]\text{,}
\end{equation}
where
$$f_{D}(\mathbf{x})=\left(1-\frac{1}{\beta}\right)\pi_{D}(\mathbf{x})+\frac{1}{\beta}\mathbf{x}\text,$$
$$f_{C}(\mathbf{x})=\left(1+\frac{1}{\beta}\right)\pi_{C}(\mathbf{x})-\frac{1}{\beta}\mathbf{x}\text.$$
In this paper we use only $\beta=1$.
A difference map search proceeds by starting from a random initial configuration $\mathbf{x}_0$
and iteratively applying the difference map: $\mathbf{x}_i=\operatorname{DM}(\mathbf{x}_{i-1})$
\cite{DM-PNAS}.
When the map reaches a fixed point $\mathbf{x}_{fp}$, a solution is obtained by
\begin{equation}
\mathbf{x}_{sol}=\pi_C\left( f_D(\mathbf{x}_{fp})\right)=\pi_D\left( f_C(\mathbf{x}_{fp})\right)\text{.}
\end{equation}
Notice that the ability to obtain a solution from any fixed point of the map, due to the cancelation
of the two bracketed terms in \eqref{DMeq}, relies on the definition of the problem in terms of
only two simple constraints.
For a given iterate $\mathbf{x}_i$, the terms $\pi_{C}\left( f_{D}(\mathbf{x}_i)\right)$ and
$\pi_{D}\left( f_{C}(\mathbf{x}_i)\right)$ provide two estimates of the solution,
each satisfying one of the two constraints. We call these respectively the $C$- and $D$-estimates
of the solution at the $i$th iteration. The distance between the two estimates
is the error $\epsilon$ and the search terminates when the error converges to zero.

To summarize, a simple difference map solver for continuous constraints
would consist of the following simple steps:
\begin{enumerate}
\item
Initialize the iterate $\mathbf{x}$ to a random configuration.
\item
Compute the two estimates of the solution $\mathbf{x}_C\leftarrow\pi_{C}\left( f_{D}(\mathbf{x})\right)$
and $\mathbf{x}_D\leftarrow\pi_{D}\left( f_{C}(\mathbf{x})\right)$.
\item
Compute the error $\epsilon\leftarrow||\mathbf{x}_C-\mathbf{x}_D||$. If it is below a predefined
convergence threshold, the search terminates, and the solution is given by
$\mathbf{x}_C\approx\mathbf{x}_D$.
\item
Advance the iterate $\mathbf{x}\leftarrow\mathbf{x}+\beta(\mathbf{x}_D-\mathbf{x}_C)$. Start the
next iteration at Step 2.
\end{enumerate}

\section{\label{sec2}Constraints}

\subsection{\label{sec22}Periodic sphere packing and kissing}

A periodic packing of equal-sized spheres (radius $r$) in $d$ dimensions can be generated by
the action of a lattice $\Lambda$ on a set of $p$ primitive spheres.
Let $P$ be the set of centers of the primitive spheres.
We define a \textit{generating matrix} of the packing as a $(d+p)\times d$ matrix $\mathrm{M}$
whose first $d$ rows are a
set of generators of $\Lambda$ and whose remaining $p$ rows are the vectors in the set $P$.
Combining these quite different sets of configuration variables into a single matrix serves to remind
us that at the highest level of our search algorithm both sets are treated in a uniform manner by
the projection operators. The detailed constraints, of course, distinguish among the two parts
of $\mathrm{M}$, which we denote $\mathrm{M}_0$ (lattice generators) and $\mathrm{M}_1$
(primitive sphere centers). The set of all the centers of spheres in the
packing is then the Minkowski sum
\begin{align}
\Lambda+P & =\{\mathbf{b}_0\mathrm{M}_0+\mathbf{y} \colon \mathbf{b}_0\in\mathbb{Z}^d,~\mathbf{y}\in P\}\\
&=\{\mathbf{b}\mathrm{M} \colon \mathbf{b} \in \mathbb{Z}^d\oplus E_{p} \} \text, \nonumber
\end{align}
where $E_p$ is the set of coordinate-permutations of the $p$-dimensional vector $(1,0,0,\ldots,0)$.
The space $\mathbb{R}^{(d+p)\times d}$ of generating matrices takes the role of the
physical configuration space $\Psi$.

A matrix $\mathrm{M}$ generates a valid packing if the centers of any two sphere of the packing, 
$\mathbf{b}_1\mathrm{M}$ and $\mathbf{b}_2\mathrm{M}$, are separated 
at least by a distance of $2r$ when $\mathbf{b}_1\neq\mathbf{b}_2$.
Each choice of $\mathbf{b}_1$ and $\mathbf{b}_2$ generates a constraint on the matrix $\mathrm{M}$
\begin{equation}
||\mathbf{b}_1\mathrm{M}-\mathbf{b}_2\mathrm{M}||\ge 2r\text,
\end{equation} 
which we call an {\it exclusion constraint}. Note that there are infinitely many independent
exclusion constraints (constraints with $\mathbf{b}_1-\mathbf{b}_2=\mathbf{b}'_1-\mathbf{b}'_2$
are not independent). However, for any non-degenerate
matrix $\mathrm{M}$ only finitely many independent exclusion constraints are violated or are even
remotely close to being violated. In practice, only those constraints need be
tested in our computations. We call those constraints the \emph{relevant} exclusion
constraints (let there be $n$ of them), and we define a $2n\times (d+p)$ matrix 
$\mathrm{A}$ whose rows $\mathbf{a}_{2i-1}$ and $\mathbf{a}_{2i}$ are the vectors
$\mathbf{b}_1$ and $\mathbf{b}_2$ related to the $i$th relevant exclusion constraint.
We discuss below how the relevant constraints are identified. 

The linear map $A: \mathrm{M}\mapsto \mathrm{X}=\mathrm{A}\mathrm{M}$ is a map from the
physical configuration space $\Psi$ to a larger-dimensional space, $\Omega=\mathbb{R}^{2n\times d}$,
which we use as the formal configuration space. As before, since the map $A$ is not surjective,
only a subset (a linear subspace, in fact) of formal configurations have a corresponding
generating matrix in the physical configuration space. The choice of constraints below will
guarantee that solutions belong to this subset.
The size of the configuration space grows as the number of relevant independent exclusion
constraints, which is the number of independent near neighbor pairs for which
overlap needs to be actively avoided.
Notice that each relevant exclusion constraint can now be written in terms of a private
set of variables. Specifically, each row of the matrix $\mathrm{X}$ corresponds to the
position of one particle, and the $i$th relevant exclusion constraint is given by
\begin{equation}
D_i = \{\mathrm{X}\in \Omega \colon ||\mathbf{x}_{2i-1}-\mathbf{x}_{2i}||\ge2r \}\text.
\end{equation}
The intersection of all the relevant exclusion constraints forms our ``divide'' constraint,
\begin{equation}
D = \{\mathrm{X}\in \Omega \colon ||\mathbf{x}_{2i-1}-\mathbf{x}_{2i}||\ge2r \text{ for } i=1,\ldots n\}\text.
\end{equation}
Each set of private variables associated with one exclusion constraint 
is composed of the coordinates of replicas of two particles,
and we call these two replicas a {\it replica pair}.

As mentioned in Section \ref{sec1a1}, the confinement constraint of finite packing problems
is replaced in the case of periodic packings with a constraint on the density of the packing.
The density of a packing generated by a matrix $\mathrm{M}$ is given by the density of
the unit cell, whose volume is $|\det \mathrm{M}_0|$ and which contains $p$ particles
of volume $V_1$:
\begin{equation}
\phi = \frac{p V_1}{|\det \mathrm{M}_0|}\text.
\end{equation}
Therefore, if we wish to find a packing of density $\phi\ge\phi_{\text{target}}$, the
{\it density constraint} on the generating matrix will be
\begin{equation}
K_\text{density} = \{\mathrm{M}\in\Psi \colon |\det \mathrm{M}_0|\le V_\text{target}\}\text,
\end{equation}
where $V_\text{target} = p V_1/\phi_\text{target}$.

As in the example of Section \ref{sec1a1}, since the map $A$ is not surjective,
a general element $\mathrm{X}\in\Omega$ of the
formal configuration space does not correspond to a well-defined physical configuration.
We therefore impose a constraint that requires $\mathrm{X}$ to lie in the range of $A$.
In the context of the PDC construction we call this the {\it lattice constraint} because
it requires different periodic images of a primitive particle to lie on the points of a lattice,
and requires that lattice to be the same for all primitive particles (up to translation).
Again, as in Section \ref{sec1a1}, we combine the lattice constraint with the density constraint
to form the ``concur'' constraint:
\begin{align}
\label{Concur}
C =& A(K_\text{density})\\
=& \{\mathrm{X}=\mathrm{A}\mathrm{M}\in \Omega \colon |\det M_0| \le V_\text{target}\}\text.\nonumber
\end{align}

\begin{figure}

\includegraphics[scale=0.29]{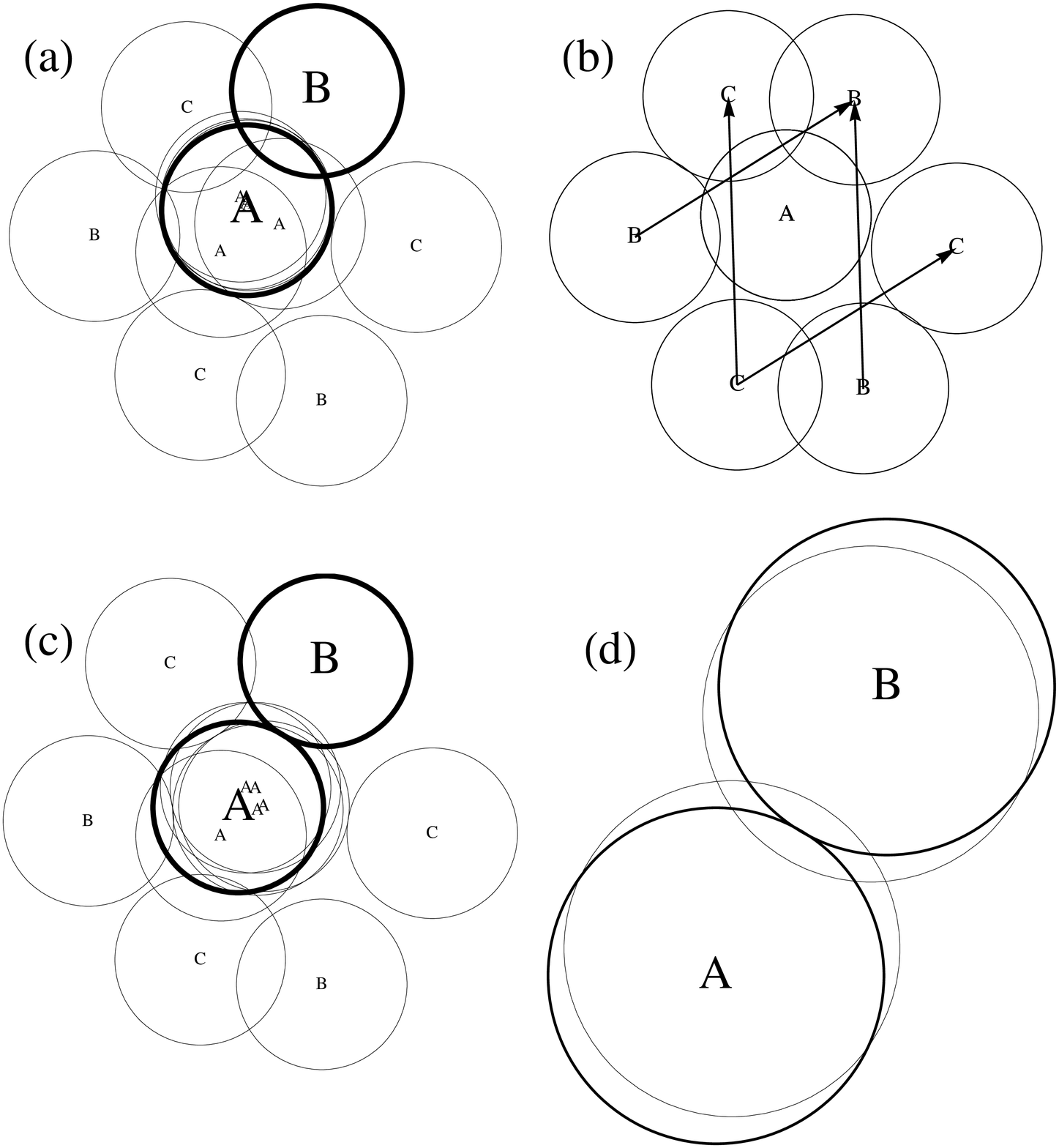}

\caption{
\label{fig1}
An illustration of the ``divide'' and ``concur'' projections in the two-dimensional sphere packing
problem with $p=3$. (a) A hypothetical configuration of six replica pairs
involving the primitive disk A. Disks with the same
letter marking  their centers are replicas of the same primitive disk (as for disk $A$) or of its 
 lattice translates (as for disks $B$ and $C$). One replica pair, violating its exclusion constraint,
 is emphasized. (b) The output of the ``concur'' projection: the closest configuration to (a) such that
 all replicas of a particular primitive disk lie on top of each other, or a lattice translation
 apart (arrows), and such that those lattice translations define a lattice with a sufficiently
 small unit cell volume. This projection is a modification of the ``concur'' projection depicted
 in Figure \ref{fig0}d,g. (c) The output of the ``divide'' projection: the closest configuration to (a) such
 that no replica pair violates its exclusion constraint. This is identical to the ``divide'' projection
 depicted in \ref{fig0}d-f.  Detail: (d) the emphasized replica pair
 before the ``divide'' projection (thin-outline disks) and after (thick-outline disks) isolated for clarity.  
 }

\end{figure}

With these definitions of the constraint sets, $\mathrm{X}=\mathrm{A}\mathrm{M}\in C\cap D$
if and only if $\mathrm{M}$ generates a periodic packing of density $\phi\ge\phi_\text{target}$. 
The action of the projections $\pi_D$ and $\pi_C$ to the two constraints is illustrated in Figure \ref{fig1}
and Sections \ref{sec31} and \ref{sec32} discuss how the projections are computed efficiently. 

The basic operations of the search --- projections --- depend directly on the metric defined on the
formal configuration space. Therefore, the choice of metric affects
both the complexity of implementing the projection and the search dynamics.
The simplest choice for the metric is the distance induced from the Frobenius (Euclidean) norm
\begin{equation}
||\mathrm{X}_1-\mathrm{X}_2||^2_F = \operatorname{trace}\left((\mathrm{X}_1-\mathrm{X}_2)(\mathrm{X}_1-\mathrm{X}_2)^T\right)\text.
\end{equation}
This choice of metric amounts to giving all replicas of a particle equal weight in influencing
its consensus position in the ``concur'' projection. We can use a slightly different Euclidean metric, given by
\begin{equation}
||\mathrm{X}_1-\mathrm{X}_2||^2_\mathrm{W} = \operatorname{trace}\left(\mathrm{W}(\mathrm{X}_1-\mathrm{X}_2)(\mathrm{X}_1-\mathrm{X}_2)^T\right)\text,
\end{equation}
where $\mathrm{W}$ is a diagonal matrix whose diagonal elements $w_i$ are the metric weights
of different replicas. Performance is greatly enhanced by adjusting 
the metric weights throughout the search to afford greater weight to
replica pairs that continually violate their constraints and smaller weight to replica pairs that are in low risk of violating their constraints \cite{D-C}.
Note that removing a constraint from the list of relevant
constraints (i.e., removing the corresponding pair of rows from $\mathrm{A}$ and $\mathrm{X}$) is
equivalent to setting the metric weight of its replicas to zero.
Therefore, in the course of the search we not only adjust
the weights $w_i$ of replica pairs, but also add and remove replica pairs. The details of how these
changes are applied systematically are given in Section \ref{sec33}.

This constraint formulation of the periodic sphere packing problem (finding a
periodic packing with density $\phi_{\text{target}}$) can be straightforwardly
modified to describe instead the periodic kissing number problem (finding a periodic packing with
average coordination number $\tau_{\text{target}}$). First, the ``divide'' constraint is modified so that
each replica pair must still be separated by a distance of at least $2r$, but at least
$p \tau_{\text{target}}$ replica pairs must be separated by a distance of exactly $2r$.
Second, the condition on the volume of the unit cell
is dropped from the ``concur'' constraint. Projections to these modified constraints are also given in
Section \ref{sec3}. 

\subsection{\label{sec23}Convex polytope packing}

The symmetry of the spherical particle allows its configuration to be described solely by the
position of its center. In the case of a general convex particle, the variables of the configuration
space need to include information also about the orientation of the particle. One possible
description of the particle assigns variables separately to the position of its centroid and
to the description of the rotation about the centroid (e.g., a rotation matrix or a quaternion).
In this paper, 
however, we find it more convenient to describe  convex polytopes by reference to
the positions of their vertices. Therefore, a polytope with $v$ vertices is represented by
a $v\times d$ vertex matrix and is given by the convex hull of these vertices. 
Although the configuration of a single particle is no
longer represented by a single vector but by a matrix composed of $v$ vectors,
it is convenient to treat these matrices as vectors, which we typeset as
bold-face upper-case
Latin letters (e.g., $\mathbf{X}$ for the vertex matrix of the polytope
$K=\operatorname{conv}\mathbf{X}=\operatorname{conv}\{\mathbf{x}_i \colon i=1,\ldots v\}$),
and to construct matrices whose rows are such vectors.
A translation by $\mathbf{t}$ of a polytope $\operatorname{conv}\mathbf{X}$ is given by
$\operatorname{conv}(\mathbf{X}+\mathbf{c}^T\mathbf{t})$, where $\mathbf{c}^T$ is a column vector
of unit elements and $\mathbf{c}^T\mathbf{t}$ is the translation matrix corresponding to the
translation vector $\mathbf{t}$.
Similarly, a rotation is given by $\operatorname{conv}(\mathbf{X}\mathrm{R})$,
where $\mathrm{R}$ is a $d\times d$ orthogonal matrix.

A periodic packing is again generated
by the action of a lattice $\Lambda$ on a set of $p$ primitive polytopes whose vertex matrices form
the set $P$.
The set of all vertex matrices of polytopes in the packing is the Minkowski sum
\begin{align}
\mathbf{c}^T\Lambda+P&=\{\mathbf{b}_0\mathrm{M_0}+\mathbf{Y} \colon \mathbf{b}_0\in\mathbb{Z}^d,~\mathbf{Y}\in P\}\\
&=\{\mathbf{b}\mathrm{M} \colon \mathbf{b} \in \mathbb{Z}^d\oplus E_{p} \}\text,\nonumber
\end{align}
where $\mathrm{M}$ is a generating matrix of the packing, whose first $d$ rows (comprising 
$\mathrm{M}_0$) are translation matrices generating $\Lambda$, and whose
remaining $p$ rows (comprising $\mathrm{M}_1$) are the vertex matrices of the set $P$.
The space of generating matrices $\Psi=\mathbb{R}^{(d+p)\times(v\times d)}$ is the
physical configuration space.

Each exclusion constraint between two particles of the packing requires the convex hulls
$\operatorname{conv}(\mathbf{b}_1\mathrm{M})$ and
$\operatorname{conv}(\mathbf{b}_2\mathrm{M})$ not to overlap for any
$\mathbf{b}_1\neq\mathbf{b}_2$. To construct the formal configuration space
we again form one replica pair for the particles involved in each relevant
exclusion constraint, which gives $\Omega=\mathbb{R}^{2n\times(v\times d)}$.
The map $A$ from physical configurations to formal configurations is given
by the matrix $\mathrm{A}$ whose rows  $\mathbf{a}_{2i-1}$ and $\mathbf{a}_{2i}$ are the vectors
$\mathbf{b}_1$ and $\mathbf{b}_2$ related to the $i$th relevant exclusion constraint.
The ``divide'' constraint is given by the intersection of all the relevant exclusion
constraints, each expressed in terms of its private replica pair:
\begin{align}\label{dividepolytopes}
D = \{\mathrm{X}\in \Omega \colon \operatorname{int}(\operatorname{conv}\mathbf{X}_{2i-1}\cap\operatorname{conv}\mathbf{X}_{2i})=\emptyset & \\
\text{for } i=1,2,\ldots n&\}\text.\nonumber
\end{align}

In addition to the lattice constraint and the density constraints, which combine in Section \ref{sec22}
to form the ``concur'' constraint, in the case at hand we must include a third constraint, the {\it
rigidity constraint}. The primitive particles of a packing generated by a general matrix $\mathrm{M}$
are only constrained in their number of vertices, not in the arrangement of those vertices. However,
we are interested only in packing where all the particles are congruent with a given shape,
and so we impose the constraint on $\mathrm{M}$ that the vertices of its primitive particles
are obtained from the vertices of the given particle by a rigid motion:
\begin{align}
K_\text{rigidity} = \{\mathrm{M}\in \Psi \colon \mathbf{Y} = \mathbf{Y}^{(0)}\mathrm{R}_i + \mathbf{c}^T\mathbf{t}_i &\\
\text{ for all $p$ rows $\mathbf{Y}$
of $\mathrm{M}_1$}&\}\text{,}\nonumber
\end{align}
where $\mathbf{Y}^{(0)}$ is the vertex matrix of the given particle. The ``concur'' constraint
$C=A(K_\text{density}\cap K_\text{rigidity})$ is given by combining the density and rigidity
constraints on the generating matrix with the lattice constraint.
The result of constructing the ``divide'' and ``concur'' constraints is that a formal configuration
satisfies both of them if and only if it corresponds to a generating matrix 
in the physical configuration space which yields a packing of the given particle with the desired
density. 
Table \ref{tab1} summarizes the $D-C$ constraints for the three problems discussed and Section \ref{sec3} describes in detail the projections to these constraints.

\begin{table}
\centering

\begin{tabular}{|l|l|l|}
\hline
& ``divide'' constraint & ``concur'' constraint \\
\hline
sphere & $||\mathbf{x}_{2i-1}-\mathbf{x}_{2i}||\ge2r$ for&$\mathrm{X}=\mathrm{A}\mathrm{M}$  \\
packing & all $n$ replica pairs& $\mathrm{M}\in K_\text{density}$\\
\hline
kissing & $||\mathbf{x}_{2i-1}-\mathbf{x}_{2i}||\ge2r$ for& $\mathrm{X}=\mathrm{A}\mathrm{M}$\\
number & all $n$ replica pairs and & \\
&$=2r$ for  $p \tau_\text{target}$ pairs & \\
\hline
polytope & convex hulls of $\mathbf{X}_{2i-1}$ & $\mathrm{X}=\mathrm{A}\mathrm{M}$\\
packing &and $\mathbf{X}_{2i+1}$ non-overlap-& $\mathrm{M}\in K_\text{density}$\\ 
 &ping for all $n$ replica & and $\mathrm{M}\in K_\text{rigidity}$ \\
& pairs & \\
\hline

\end{tabular}

\caption{\label{tab1}A summary of the $D-C$ constraints for periodic sphere packing, the average kissing number
problem, and polytope packing. The ``divide'' constraint encompasses the relevant exclusion constraints, while the ``concur'' constraint encompasses, where applicable, the density, rigidity, and
lattice constraints.}

\end{table}

\section{\label{sec3}Implementation}

\subsection{\label{sec31}``Divide'' projections}

\subsubsection{\label{sec311}Sphere packing and kissing}

In order to implement an iterated difference map search, whose iterations are given by \eqref{DMeq},
we must implement efficient projections to the ``divide'' and ``concur'' constraints. These
implementations are the subject of Sections \ref{sec31} and \ref{sec32}. In the course of the
search, considerations of efficiency require certain changes to the formal configuration space --
specifically adding and removing replica pairs, changing metric weights, and lattice reduction.
In section \ref{sec33} we discuss when and how these changes are applied. 

In the case of sphere packing, the ``divide'' constraint simply requires
that the centers of the two spheres comprising each replica pair be a certain distance apart.
This is obtained by applying equation \eqref{sphereproj} to each replica pair.
Note that the ``divide'' projection acts independently on each replica pair, and 
since the metric weight of all variables specifying one replica pair are equal,
the metric weights have no influence on this projection.
The action of this projection is illustrated in Figure \ref{fig1}.
For the kissing number problem, the first case of \eqref{sphereproj}
is also used if the replica pair is one of the $p \tau_{\text{target}}$ closest replica pairs.

\subsubsection{\label{sec312}Convex polytope packing}

\begin{figure}

\includegraphics[scale=0.41]{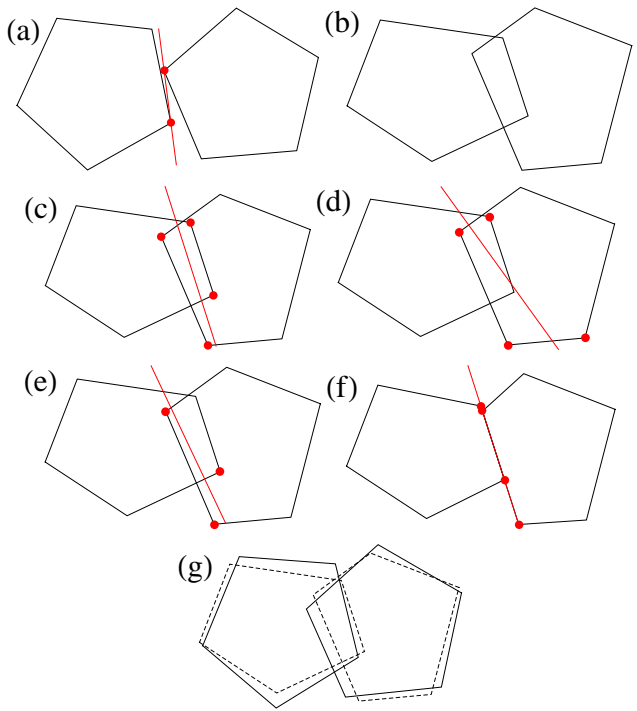}

\caption{\label{fig2}
An illustration of the polytope exclusion and rigidity constraint projections for the case of
regular pentagons. (a) The pentagons are non-overlapping, as demonstrated by the existence of
a separating axis that passes through one vertex of each pentagon (red dots).
(b) A hypothetical situation of overlapping pentagons. No axis, and particularly
no axis passing through one vertex of each pentagon, separates the two sets of vertices.
(c) For the input pentagons in (b), 
the subset $S=T$ (red dots) that
minimizes $\delta(S)^2$, the sum of squared-distances to the least-squares axis (red line) while
satisfying that the latter separates the remaining vertices.
(d) Another choice of $S$ that yields a valid separating axis, but a larger
sum of squared-distances. (e) A choice of $S$ that yields an axis that fails
to separate the remaining vertices. (f) Using
$T$ found in (c), the output of the projection to the exclusion constraint is determined by moving
the points of $T$ onto their least-squares axis. (g) The output (solid line) of the rigidity
projection for the input pentagons (dashed) in (b). 
 }
\end{figure}

Although identifying and resolving overlaps between two spheres is straightforward,
the same task is more challenging in the case of other convex objects, where more degrees of
freedom come into play. The literature on the topic of
detecting overlaps (collisions) between polyhedral solids is extensive
(driven in part by applications in
computer graphics), and many efficient techniques exist for checking
whether two convex polyhedra, $\operatorname{conv}X_1$ and
$\operatorname{conv}X_2$, overlap (see e.g., \cite{collision1,collision2}).
In our case, as we are interested in computing the projection to the exclusion constraint,
we also need to determine the distance-minimizing resolution of the overlap. That is,
we must find the smallest displacement of the vertices such that the new polyhedra 
do not overlap.
As far as we have been able to determine, there is not an established, efficient computational
method developed for this specific problem. The method we provide here is efficient enough for the
purpose of packing polyhedra with a small number of vertices, but the
computation time required grows exponentially with the number of vertices. A more efficient resolution
method for particles with more vertices and for smooth particles is currently in development.

The method relies on the separating plane theorem: the convex hulls of two
sets of vertices in $\mathbb{R}^d$ do not overlap
if and only if there is a $(d-1)$-dimensional plane that separates the two sets, so that each is
contained in a different half-space. The theorem can be made even stronger
by specifying that the separating plane can always be chosen to contain $d$ vertices from the
given sets, including at least one from each. Therefore, one can check whether two 
polytopes overlap by checking whether they are separated by any of the planes
defined by any such subset of vertices (Figure \ref{fig2}a--b). If the polytopes
are non-overlapping, the resolution leaves the vertices unchanged. 

If the polytopes overlap, we must find the smallest displacement of their
vertices that resolves the overlap. In the resolved configuration there is a separating
plane, $V_{sp}=\{\mathbf{r}\in\mathbb{R}^d \colon \hat{\mathbf{n}}_{sp}\cdot\mathbf{r}=h_{sp}\}$,
that separates the two sets of vertices. As a consequence of distance minimization,
the only vertices moved in the course of the resolution are the vertices which lie on $V_{sp}$
in the resolved configuration. Let $T$ and $T'$ be the pre- and post-resolution positions, respectively,
of those vertices that are displaced during the resolution. Therefore, $T'$ is the set of points in
$V_{sp}$ closest to the points of $T$:
\begin{equation}
T'=\{ \mathbf{r} + (h_{sp}-\mathbf{r}\cdot\hat{\mathbf{n}}_{sp})\hat{\mathbf{n}}_{sp} \colon \mathbf{r}\in T\}\subseteq V_{sp}\text{.}
\end{equation}
The squared norm of the resolution displacement is 
\begin{equation}
\sum_{\mathbf{r}\in T} (h_{sp}-\mathbf{r}\cdot\hat{\mathbf{n}}_{sp})^2\text{.}
\end{equation}

For any set of points $S$ there is at least one plane
$V=\{\mathbf{r}\in\mathbb{R}^d \colon \hat{\mathbf{n}}\cdot\mathbf{r}=h\}$
that minimizes the sum of squared distances
\begin{equation}\label{planedist}
\sum_{\mathbf{r}\in S} (\hat{\mathbf{n}}\cdot\mathbf{r}-h)^2\text.
\end{equation}
We call such a plane a least-squares plane of $S$.
The separating plane of the resolved configuration is always a least-squares plane of $T$.
If this were not the case, a small tilting of the separating
plane towards such a least-squares plane
(with a corresponding movement of the points in $T'$) would result in a resolution by
a smaller displacement. In order to resolve an overlap between polytopes $\operatorname{conv}X_1$
and $\operatorname{conv}X_2$, we therefore have to solve a discrete problem: 
among all subsets $S$ of $X_1\cup X_2$ with a least-squares
plane separating the remaining vertices $X_1\setminus S$ from $X_2\setminus S$, find the one with the minimal sum of squared distances \eqref{planedist}. This is the set $T$ (Figure \ref{fig2}c--f).

The least-squares plane $V_{ls}$ of a set $S$ is
determined by minimizing the sum of squared distances \eqref{planedist}.
Note that for a fixed normal direction $\hat{\mathbf{n}}$, the value of
$h$ that minimizes the sum is $h=\hat{\mathbf{n}}\cdot\overline{\mathbf{r}}$, where
$\overline{\mathbf{r}} = \sum_{\mathbf{r}\in S}\mathbf{r}/|S|$ is the centroid of $S$. Therefore, we wish to minimize
\begin{equation}
\sum_{\mathbf{r}\in S} [\hat{\mathbf{n}}\cdot(\mathbf{r}-\overline{\mathbf{r}})]^2 = 
\hat{\mathbf{n}} \left[ \sum_{\mathbf{r}\in S} (\mathbf{r}-\overline{\mathbf{r}})^T(\mathbf{r}-\overline{\mathbf{r}}) \right ] \hat{\mathbf{n}}^T\text,
\end{equation}
the minimum of which is equal to the smallest eigenvalue of the symmetric matrix
$\sum_{\mathbf{r}\in S} (\mathbf{r}-\overline{\mathbf{r}})^T(\mathbf{r}-\overline{\mathbf{r}})$.
The minimum is realized when $\hat{\mathbf{n}}$ is the corresponding eigenvector.
Degenerate cases with equal lowest eigenvalues occur, but they do not pose a problem:
whenever an optimal separating plane occurs as a degenerate least-squares
plane of some set $S$, its degeneracy implies that there is a least-squares plane of $S$
which also includes an extra vertex; this plane will be equally optimal and will occur
as a less degenerate least-squares plane of a superset $S'\supseteq S$. Therefore,
the optimal least-squares plane always occurs as a non-degenerate least-squares
plane of a set $S$.

To summarize, the overlap detection and resolution algorithm
consists of three steps (illustrated in Figure \ref{fig2}a--f):
\begin{enumerate}
\item Consider all subsets $S\subseteq X_1 \cup X_2$ of size $|S|=d$ with
at least one point from each polytope.
Let $V=\{\mathbf{r}\in\mathbb{R}^d \colon \hat{\mathbf{n}}\cdot\mathbf{r}=h\}$ be
a plane that includes $S$. For each $S$ let
\begin{align}
\Delta_+^2(S) & = 
\sum_{\substack{\mathbf{x}\in X_1 \\ \hat{\mathbf{n}}\cdot\mathbf{x}>h}}
(\hat{\mathbf{n}}\cdot\mathbf{x}-h)^2 + 
\sum_{\substack{\mathbf{x}\in X_2 \\ \hat{\mathbf{n}}\cdot\mathbf{x} \le h}}
(\hat{\mathbf{n}}\cdot\mathbf{x}-h)^2 \text{,} \\
\Delta_-^2(S) & = 
\sum_{\substack{\mathbf{x}\in X_1 \\ \hat{\mathbf{n}}\cdot\mathbf{x} \le h}}
(\hat{\mathbf{n}}\cdot\mathbf{x}-h)^2 + 
\sum_{\substack{\mathbf{x}\in X_2 \\ \hat{\mathbf{n}}\cdot\mathbf{x} > h}}
(\hat{\mathbf{n}}\cdot\mathbf{x}-h)^2 \text{,} \\
&\Delta^2(S)  = \operatorname{min}(\Delta_+^2(S),\Delta_-^2(S)) \text{,}
\end{align}
and let
\begin{equation}
\label{Delta}\Delta^2  = \underset{S}{\operatorname{min}} \Delta^2(S) \text{.}
\end{equation}
$\Delta^2$ provides a measure for the interpenetration of the two polytopes.
If $\Delta^2=0$, then a separating plane exists, 
the input polytopes do not overlap, and the algorithm
ends here by returning the original
vertex positions $X_1$ and $X_2$. If $\Delta^2>0$, 
the polytopes overlap and the algorithm continues to Step 2.
\item Consider all subsets $S\subseteq X_1 \cup X_2$ of size $|S|>d$ with
at least one point from each polytope.
Let $V=\{\mathbf{r}\in\mathbb{R}^d \colon \hat{\mathbf{n}}\cdot\mathbf{r}=h\}$
be a least-squares plane of $S$.
If the plane separates the vertex sets with the points of $S$ removed ---   
$X_1\setminus S$ and $X_2\setminus S$ --- 
let $\delta^2(S)$ be the sum of squared-distances from $S$ to the plane. Otherwise, let $\delta^2(S)=\infty$.
Among the subsets $S$ considered, let $T$ be the subset that minimizes $\delta^2(S)$ and 
$V_T=\{\mathbf{r}\in\mathbb{R}^d \colon \hat{\mathbf{n}}_T\cdot\mathbf{r}=h_T\}$ be its
associated least-squares plane. As $\delta^2(S)<\infty$ if $S$ contains all vertices, the minimum
is always finite. Continue to Step 3.
\item The sets of vertices returned are given by $X_1'$ and $X_2'$, wherein
$\mathbf{x}'\in X_1'\cup X_2'$ is given by
\begin{equation}
\mathbf{x}'=\begin{cases} \mathbf{x} & \text{if $\mathbf{x}\not\in T$ }\\
\mathbf{x} + (h_T - \mathbf{x}\cdot\hat{\mathbf{n}}_T)\hat{\mathbf{n}}_T & \text{if $\mathbf{x}\in T$, }
\end{cases}
\end{equation}
where $\mathbf{x}\in X_1\cup X_2$ is the corresponding original vertex position.

\end{enumerate}

The projection $\pi_D(\mathrm{X})$ to the ``divide'' constraint \eqref{dividepolytopes},
of an input matrix $\mathrm{X}$
comprised of pairs of vertex matrices $\mathbf{X}_{2i-1}$ and $\mathbf{X}_{2i}$, is then
achieved by applying the above algorithm independently to all $i=1,\ldots n$ pairs.

\subsection{\label{sec32}``Concur'' projections}

\subsubsection{\label{sec321}Lattice constraint}

All the ``concur'' constraint sets described in this paper are of the form
\begin{equation}
\mathrm{C}=A(K)=\{\mathrm{X}=\mathrm{A}\mathrm{M}\in\Omega \colon \mathrm{M}\in K\}
\end{equation}
where $\mathrm{A}$ is constant, and $\mathrm{M}$ is variable,
but must satisfy a constraint $\mathrm{M}\in K$. The projection then is given by
\begin{equation}
\pi_C:\mathrm{X}\mapsto\mathrm{X'} = \mathrm{A}\mathrm{M}\text,
\end{equation}
where $\mathrm{M}$ realizes the minimum over $K$ of the distance
\begin{equation}\label{concdist}
||X-X'||^2 = \operatorname{trace}\left(\mathrm{W}(\mathrm{X}-\mathrm{A}\mathrm{M})(\mathrm{X}-\mathrm{A}\mathrm{M})^T\right)\text.
\end{equation}

Absent any constraints on $\mathrm{M}$ (as for example in the
``concur'' constraint for the kissing number problem,
where $K=\Psi$),
the solution would be given by
\begin{equation}
\overline{\mathrm{M}}=(\mathrm{A}^T\mathrm{W}\mathrm{A})^{-1}
\mathrm{A}^T\mathrm{W}\mathrm{X}\text.
\end{equation}
This can easily be seen by writing $\mathrm{M}=\overline{\mathrm{M}}+\delta\mathrm{M}$, which gives
\begin{align}
||X-X'||^2 &= \mathrm{trace}\left(\mathrm{W}(\mathrm{X}-\mathrm{A}\mathrm{M})(\mathrm{X}-\mathrm{A}\mathrm{M})^T\right) \nonumber \\
 &= c+\mathrm{trace}(\mathrm{W}\mathrm{A}\,\delta\mathrm{M}\,\delta\mathrm{M}^T\,\mathrm{A}^T)\nonumber\\
\label{costX} &= c+\mathrm{trace}(\mathrm{W}' \delta\mathrm{M}\,\delta\mathrm{M}^T)\text,
\end{align}
where $\mathrm{W}' =\mathrm{A}^T\mathrm{W}\mathrm{A}$ and
the constant term $c$ does not depend on $\delta\mathrm{M}$. The second term is non-negative, and 
when $\mathrm{M}$ is unconstrained, \eqref{costX} is minimized by letting
$\mathrm{M}=\overline{\mathrm{M}}$. Additionally, we have just reduced the 
constrained case to the problem of finding $\mathrm{M}\in K$ that minimizes the cost function
\begin{equation}\label{fM}
f(\mathrm{M})=\operatorname{trace}\left(\mathrm{W}' (\mathrm{M}-\overline{\mathrm{M}})(\mathrm{M}-\overline{\mathrm{M}})^T\right)\text.
\end{equation}

This projection strategy parallels the two-step strategy used in Section \ref{sec1a1}.
First, the formal configuration $\mathrm{X}$ is projected to the range $A(\Psi)$ of
the physical configuration space, giving $\mathrm{A}\overline{\mathrm{M}}$.
Then, the projection of $\mathrm{\overline{M}}$ to the
additional constraint $K$ is performed in the physical configuration
space using the metric induced on its image in the formal configuration space.
Below, we solve the second step of this projection problem for various constraints $K$.

\subsubsection{\label{sec322}Density constraint}

In the ``concur'' constraint for the sphere packing problem, the only
constraint on the generating matrix is the density constraint.
The set of generating matrices $\mathrm{M}$ satisfying the density constraint is
\begin{equation}
K_{\text{density}}=\{M \colon |\det\mathrm{M}_0| \le V_\text{target}\}\text,
\end{equation}
where $\mathrm{M}_0$ is the generating matrix of the lattice
and is given by the first $d$ rows of $\mathrm{M}$.
If $|\det\overline{\mathrm{M}}_0| \le V_\text{target}$, then the projection to the constraint (the choice of
$\mathrm{M}$ that minimizes the cost function \eqref{fM}) is trivially $\mathrm{M}=\overline{\mathrm{M}}$.
Otherwise, since $\mathrm{M}_1$ is unconstrained, we can minimize \eqref{fM} with respect to $\mathrm{M}_1$
for a given $\mathrm{M}_0$. This yields $\mathrm{M}_1 = \overline{\mathrm{M}}_1 
- \mathrm{W'}_{11}^{-1} \mathrm{W}'_{10} (\mathrm{M}_0-\overline{\mathrm{M}}_0)$, 
where $\mathrm{W}'_{IJ}$ are the block-elements of $\mathrm{W}'$ acting on $\mathrm{M}_I$ to
the left and on $\mathrm{M}_J$ to the right.
Thus, the cost function for $\mathrm{M}_0$ is simply
\begin{equation}
f(M_0)=\operatorname{trace}\left(\mathrm{W}''(\mathrm{M}_0-\overline{\mathrm{M}}_0)(\mathrm{M}_0-\overline{\mathrm{M}}_0)^T\right)\text,
\end{equation}
where $\mathrm{W}''=\mathrm{W}'_{00}-\mathrm{W}'_{01} \mathrm{W'}_{11}^{-1} \mathrm{W}'_{10}$.

The projection becomes easier to analyze in terms of the matrix
$\mathrm{L}=(\mathrm{W}'')^{1/2}\mathrm{M}_0$. The cost function then takes
the form of the simple Frobenius distance
\begin{equation}
f(\mathrm{L})=\mathrm{trace}\left((\mathrm{L}-\overline{\mathrm{L}})(\mathrm{L}-\overline{\mathrm{L}})^T\right)\text,
\end{equation}
and the density constraint is still in the form
\begin{equation}\label{Vprime}
|\det\mathrm{L}|\le V'_\text{target}\text,
\end{equation}
where $V'_\text{target} = V_\text{target} / |\det W''|^{1/2}$.
Since the absolute value of the determinant of $\mathrm{L}$ is given by the product
of its singular values, the solution to this minimization problem is given by a matrix
$\mathrm{L}=\mathrm{U}\mathrm{\Sigma}\mathrm{V}$ with the same (right and left) singular vectors
as the matrix $\overline{\mathrm{L}}=\mathrm{U}\overline{\mathrm{\Sigma}}\mathrm{V}$, but
different singular values. 
The cost function expressed in terms of the singular values $\sigma_i$ and $\overline{\sigma_i}$ of,
respectively, $\mathrm{L}$ and $\overline{\mathrm{L}}$ takes the form
\begin{equation}
f(\mathrm{\Sigma}) = \sum_{i=1}^{d} (\sigma_i-\overline{\sigma_i})^2\text.
\end{equation}
We numerically minimize this quadratic function subject to the density constraint \eqref{Vprime}.
Through back substitution we then have the matrix $\mathrm{M}$ that minimizes \eqref{concdist} and 
$\pi_C(X)=\mathrm{X}'=\mathrm{A}\mathrm{M}$. 

\subsubsection{\label{sec323}Rigidity constraint}

In the ``concur'' constraint for the polytope packing problem,
an additional constraint on the generating matrix $\mathrm{M}$ 
is that the primitive polytopes that make up $\mathrm{M}_1$ are congruent
with a given polytope. The generating matrix is then constrained
to the set
\begin{equation}
K = K_{\text{density}}\cap K_{\text{rigidity}}
\end{equation}
where
\begin{align*}
K_{\text{density}} = \{& \mathrm{M} \colon |\det \mathrm{M}_0| \le V_\text{target}\}\text, \\
K_{\text{rigidity}} = \{& \mathrm{M} \colon \mathbf{Y} = \mathbf{Y}^{(0)}\mathrm{R}_i + \mathbf{c}^T\mathbf{t}_i
\\ &\text{ for all $p$ rows $\mathbf{Y}$
of }\mathrm{M}_1\}\text{.}
\end{align*}
To calculate the projection $\pi_C(\mathrm{X})$, the cost function \eqref{fM} must be minimized over $K$.
However, since the off-diagonal block $\mathrm{W}'_{01}$
couples the lattice parameters $\mathrm{M}_0$ to the primitive particle parameters $\mathrm{M}_1$,
this minimization is complicated. Instead of exact minimization,
we employ a two-step heuristic method, which results in an approximate projection.

In the first step, we calculate the matrix $\mathrm{M}'\in K_\text{density}$
that minimizes the cost function, as in Section
\ref{sec322}. Then, in the second step, we calculate the matrix $\mathrm{M}\in K$ by applying to
each row $\mathbf{Y}$ of $\mathrm{M}'_1$
the smallest change so that it becomes a vertex matrix of a polytope congruent with the reference polytope.
The second step is achieved by finding the rigid motion applied to the reference polytope which brings
its vertices as close as possible to the vertices of $\mathbf{Y}$ as measured by the sum of squared distances
(Figure \ref{fig2}g). 
The problem of finding the rigid motion that brings one given list of points closest to another
given list, sometimes known as the problem of absolute orientation, occurs frequently in a variety
of fields (e.g., in calculating RMSD between two conformations of a biomolecule)
and several efficient methods for its solution have been developed (see \cite{RMSDq,RMSDm}).

The output of the approximate projection is then given by $\mathrm{X}'=\tilde\pi_C(\mathrm{X})=\mathrm{A}\mathrm{M}\approx\pi_C(\mathrm{X})$.
As $\mathrm{X}'\in C$, the approximate projection gives a configuration in the constraint set, but might not give
the closest one to the input configuration. We justify the use of the approximate projection by noting that
it is an exact projection if the off-diagonal block $\mathrm{W}'_{01}$ is zero. A non-zero off-diagonal block
is the result of correlations in the relevant exclusion constraint vectors $\mathbf{b}$ 
between the coefficients of lattice translations and the coefficients of primitive particle vertex positions.
We expect these coefficients to give uncorrelated contributions and to add up to small off-diagonal elements
due to random cancellations.
Indeed, we find that the off-diagonal block is small in comparison with the diagonal blocks,
and we expect our heuristic to yield a good approximate projection.

\subsection{\label{sec33}Formal configuration space maintenance}

In our discussion of the choice of metric in Section \ref{sec2}, we discussed the ideas of
dynamically readjusting the metric (through the weights $w_i$ of the various replicas) and of
removing and adding replicas (removing replicas is formally equivalent to setting their weight
to zero). The latter is necessary for
implementation reasons: there are infinitely many independent exclusion constraints (and
therefore replicas), but we can only represent a finite number of replicas in our implementation. 
As the set of relevant constraints changes over the course of the search,
we must remove and add replicas. Our
criterion for which replicas to represent is based on the difference map's current ``concur'' estimate:
we include a replica pair for each pair of particles whose centroids in the ``concur'' estimate 
are closer than some cut-off distance. Using the generating matrix obtained in the ``concur'' projection
we can easily find all such pairs using the method of Agrell \textit{et al.} \cite{LattPt}.
The cut-off distance is chosen
so that at least all replicas that might be in risk of overlap are represented.

The problem of implementation is not the only reason we wish to limit the number of replicas we
represent. A proliferation of unnecessary replicas has the adverse effect of attenuating
the information obtained from the ``concur'' projection by diluting the influence of more
critical replicas. We observe that such replica proliferation could result not only in a
slower search, but also
in an increased tendency to become trapped in local optima.
Limiting the number of replicas is one way to avoid this effect,
but we find it useful to further amplify the information from
critical constraints by giving them greater weights \cite{D-C}.
We perform the weight adjustments adiabatically,
that is, slowly over the course of many iterations,
by updating the weights of each replica pair according to the rule
\begin{equation}
w_i \to \frac{\tau w_i + w_i'(\mathrm{X}_c)}{\tau+1}\text,
\end{equation}
where $w_i'(\mathrm{X}_c)$ is a function that assigns replicas weights based on their configuration
in the ``concur'' estimate, and $\tau$ is a relaxation time for the replica weights in units of iterations.

In the sphere packing problem (in $d$ dimensions, with unit spheres),
we choose the weight function to be
\begin{equation}
w_i'(\mathrm{X}_c)=\begin{cases} e^{\alpha(4-||\mathbf{x}_i||^2)} & \text{if $||\mathbf{x}_{i}|| \le 2$} \\
(||\mathbf{x}_{i}||^2 - 3)^{-2-d/2} & \text{if $||\mathbf{x}_{i}|| > 2$,}
\end{cases}
\end{equation}
with $\alpha \approx 20$. The dimensional dependence is chosen so that under
the assumption of uniform density, the total weight from replicas over a certain distance
follows a dimension-independent power law. In the polytope packing problem, we similarly use
\begin{equation}\label{polyweight}
w_i'(\mathrm{X}_c)=\begin{cases} e^{\alpha \Delta_i^2} & \text{if the polytopes overlap}\\
(1+r_{i}^2-4r_{in}^2)^{-2} & \text{if not,}
\end{cases}
\end{equation}
with $\alpha \approx 10$, where $r_{in}$ is the inradius of the polytope, $r_i$ is the
centroid-centroid distance of the polytopes, and $\Delta_i^2$ is the measure
of the overlap between the polytopes defined in
\eqref{Delta}.

In addition to the maintenance of replicas, which is performed after every iteration of the difference map,
we also periodically perform a lattice reduction using the LLL algorithm \cite{LLL}.
The lattice generated by $\mathrm{M}_0$ is re-represented
using the LLL-reduced generating matrix
$\mathrm{M}_0'=\mathrm{G}_0\mathrm{M}_0$, where $\mathrm{G}_0$ is a unimodular integer
matrix. Additionally, all primitive particles whose centroids
are outside of the unit cell given by $\{\sum_i \lambda_i \mathbf{a}_i \colon -1/2\le\lambda_i<1/2\}$
are re-represented by their lattice-translate in that cell. In summary, the new 
packing generating matrix $\mathrm{M}'$ is given by
\begin{equation}
\mathrm{M}' = \mathrm{G} \mathrm{M} =
\left( \begin{array}{ccc} \mathrm{G}_0 & 0 \\ \mathrm{G}_1 & 1 \end{array} \right)\mathrm{M}\text,
\end{equation}
where $\mathrm{G}_1$ gives the lattice translations to be applied to the primitive particles.
Since the actual positions of the particles, as represented in the matrix
$\mathrm{X}=\mathrm{A}\mathrm{M}$, should be unchanged, the lattice reduction must also
be applied to the nominally constant matrix $\mathrm{A}$ ($\mathrm{A} \to \mathrm{A}' = \mathrm{A}\mathrm{G}^{-1}$). 

\section{\label{sec4}Results}

\subsection{\label{sec41}Sphere packing}

Using the PDC scheme described in the previous sections 
we perform a \textit{de novo} search for the densest lattice
($p=1$) sphere packings
in dimensions $2$---$14$. The PDC search, starting from random initial configurations,
was able to reproduce the densest packing lattices known 
for all cases, and the results of the search are summarized in Table \ref{tab2}. For dimensions $2$---$8$ the
lattices are known to be optimal, and for dimensions $9$---$14$ these results are,
to our knowledge, the first numerical
evidence from a \textit{de novo} search that the known lattices are optimal.

Note that the number of replicas is determined by the number of near neighbors of each sphere,
which rises rapidly with the number of dimensions.
This rise causes an increased computational storage cost per physical degree of freedom
in a PDC search, compared to a constant storage cost per physical degree of freedom in
a method involving a local search in the physical configuration space. However,
this rise need not affect the scaling of CPU costs, since both search methods
need necessarily check a comparable number of particle pairs for possible overlaps.

In dimensions $d=10,11,13$ there are known non-lattice packings with $p=40,72,144$ respectively
that are denser than the densest known lattices \cite{SPLAG}. 
In up to 11 dimensions, we searched for non-lattice packings with as many as $p=12$ primitive
spheres, but the searches did not produce packings denser than the lattice packings. For a
density target matching the lattice density, the searches reproduced the lattice packing, suggesting
that the lattice packing in these dimensions
is the optimal packing with a small number of spheres in the unit cell.

\begin{table}
\centering
\begin{tabular}{lllllll}
$d$ & $\Lambda_{\text{densest}}$ & $\phi^{(L)}_\text{densest}$ & $\langle N_\text{iter}\rangle$ & 
$\langle n\rangle$ & $t_\text{iter}$ & success rate \\
\hline 
2 & $A_2$ & $0.90690$ & $42$ & $11$ & $0.1 ms$ & $100/100$ \\
3 & $D_3$ & $0.74047$ & $230$ & $38$ & $0.2 ms$ & $100/100$ \\
4 & $D_4$ & $0.61685$ & $191$ & $127$ & $0.4 ms$ & $100/100$ \\
5 & $D_5$ & $0.46526$ & $308$ & $323$ & $1 ms$ & $100/100$ \\
6 & $E_6$ & $0.37295$ & $173$ & $977$ & $2 ms$ & $100/100$ \\
7 & $E_7$ & $0.29530$ & $217$ & $2740$ & $5 ms$ & $96/100$ \\
8 & $E_8$ & $0.25367$ & $99$ & $8528$ & $20 ms$ & $96/100$ \\
9 & $\Lambda_9$ & $0.14577$ & $161$ & $16314$ & $30 ms$ & $85/100$ \\
10& $\Lambda_{10}$ & $0.092021$ & $394$ & $31433$ & $70 ms$ & $47/100$ \\
11& $K_{11}$ & $0.060432$ & $421$ & $68722$ & $0.3 s$ & $54/100$ \\
12& $K_{12}$ & $0.049454$ & $397$ & $204321$ & $0.9 s$ & $55/100$ \\
13& $K_{13}$ & $0.029208$ & $577$ & $430796$ & $2 s$ & $25/100$ \\
14& $\Lambda_{14}$ & $0.021624$ & $1652$ & $1007250$ & $6 s$ & $4/10$ \\
\hline 
\end{tabular}

\caption{\label{tab2}Results of PDC searches for dense lattice packing in dimensions $d=2,\ldots 14$.
For each dimension, 100 runs from random initial conditions were performed with the
density target $\phi_\text{target}=\phi^{(L)}_\text{densest}$, the density of the densest known lattice $\Lambda_\text{densest}$ \cite{SPLAG}. The runs were limited to 5000 iterations,
and the number of converged runs is quoted in the right-most column.
For dimensions $10$ and above, each run was first
allowed to converge at a density target of $0.8\phi_\text{densest}$ and then
continued with the final target. The mean number of difference map iterations in
converged runs was $\langle N_\text{iter}\rangle$, and the mean number of relevant
exclusion constraint used was $\langle n\rangle$. Each iteration took an average
runtime of $t_\text{iter}$ on a single 3 GHz CPU.
In $d=14$ only 10 runs were performed with three intermediate targets.
}

\end{table}

\subsection{\label{sec42}Kissing number}

For the kissing number problem, PDC searches were able to reproduce the
best known lattice kissing arrangements
in dimensions $2$---$11$. In dimensions $2$---$9$, the result is known to be optimal, and for
dimensions $10$ and $11$, we are not aware of previous numerical evidence for their optimality.
Table \ref{tab3} summarizes the performance of our method.

\begin{table}
\centering
\begin{tabular}{llllll}
$d$ & $\Lambda_{\text{highest}}$ & $\tau^{(L)}_\text{highest}$ & $\langle N_\text{iter}\rangle$ & 
$\langle n\rangle$ & success rate \\
\hline 
2 & $A_2$ & $6$ & $27$ & $12$ & $100/100$ \\
3 & $D_3$ & $12$ & $54$ & $40$ & $100/100$ \\
4 & $D_4$ & $24$ & $132$ & $118$ & $98/100$ \\
5 & $D_5$ & $40$ & $163$ & $331$ & $94/100$ \\
6 & $E_6$ & $72$ & $225$ & $928$ & $64/100$ \\
7 & $E_7$ & $126$ & $597$ & $2729$ & $66/100$ \\
8 & $E_8$ & $240$ & $511$ & $6988$ & $55/100$ \\
9 & $\Lambda_9$ & $272$ & $350$ & $15604$ & $63/100$ \\
10& $\Lambda_{10}$ & $336$ & $438$ & $32203$ & $28/100$ \\
11& $\Lambda_{11}$ & $438$ & $549$ & $73766$ & $10/100$ \\
\hline 
\end{tabular}

\caption{\label{tab3}Results of PDC searches for lattice packing with high kissing number
in dimensions $d=2,\ldots 11$.
For each dimension, 100 runs from random initial conditions were performed with a target
coordination $\tau_\text{target}=\tau^{(L)}_\text{highest}$, the highest coordination number
known for a lattice of that dimension, $\Lambda_\text{highest}$ \cite{SPLAG}.
The runs were limited to 5000 iterations,
and the number of converged runs is quoted in the right-most column.
The mean number of difference map iteration in
converged runs was $\langle N_\text{iter}\rangle$, and the mean number of relevant
exclusion constraints used was $\langle n\rangle$.
}

\end{table}

\subsection{\label{sec43}Polytope packing}

\begin{figure}

\includegraphics[scale=0.38]{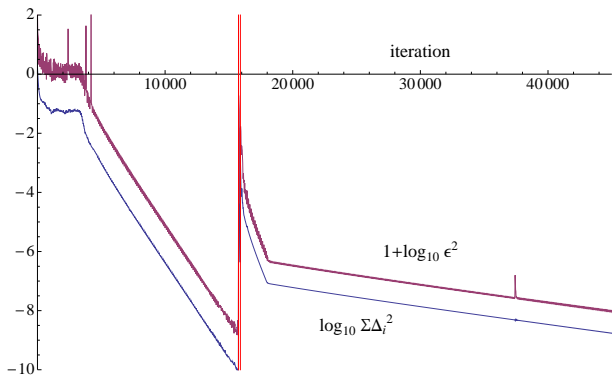}

\caption{\label{fig3}
The course of a sample run searching for dense periodic packings ($p=4$) of 
unit edge-length regular
tetrahedra, showing $\sum\Delta_i^2$, a measure of the total interpenetration between tetrahedra in the ``concur'' estimate (blue, defined in \eqref{Delta}), 
and $\epsilon^2$, the squared distance between the ``divide'' and ``concur'' estimates
(purple, shifted up for clarity), both on a logarithmic scale. The density target
for the search is started at $\phi_\text{target}=0.75$ and adjusted when the search is converged on a solution (vertical red lines) to $\phi_\text{target}=0.82$ (at iteration $15751$) and then to $\phi_\text{target}=0.8563$ (at iteration $15898$). Each iteration took 14 millisecond on average on
a single 3 GHz CPU. }
\end{figure}

By inspection of a packing of regular tetrahedra yielded by our numerical search
during early phases of its development, we were able to construct a new transitive, periodic ($p=4$)
packing of tetrahedra with a higher density ($\phi\approx0.8547$)  
than previously reported \cite{tetra}. This packing takes the form of a double lattice of
bipyramidal dimers (the union of two face-sharing tetrahedra).
The packing has since been slightly improved to a closely related, but less symmetric packing
with density $\phi\approx0.8563$ \cite{TorqUniform,ChenUniform}.
In its current form, 
our search method is able to reproduce this densest known packing reliably (fifteen
out of a hundred runs converged within the iteration limit), and Figure \ref{fig3}
shows the results of a sample run converging to this packing.

\begin{figure}

\includegraphics[scale=0.66]{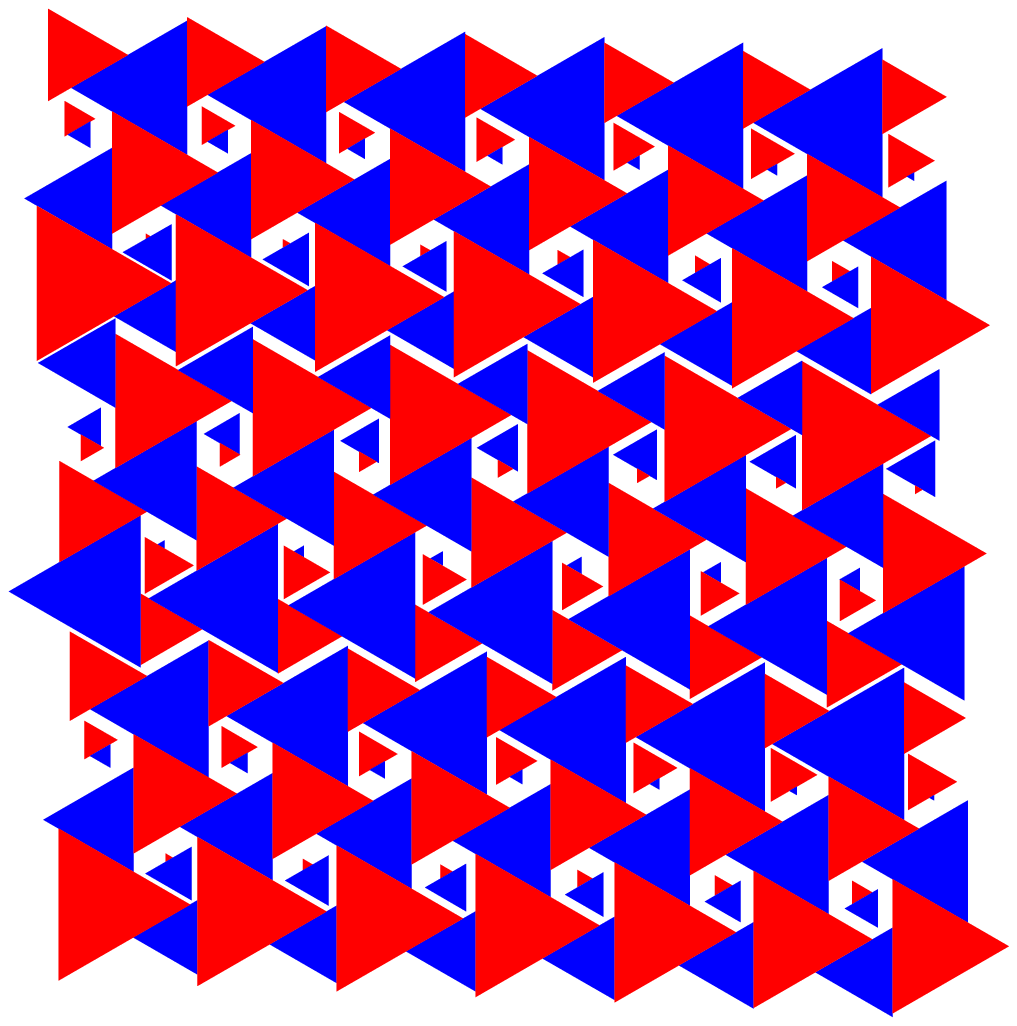} \\
\includegraphics[scale=0.55]{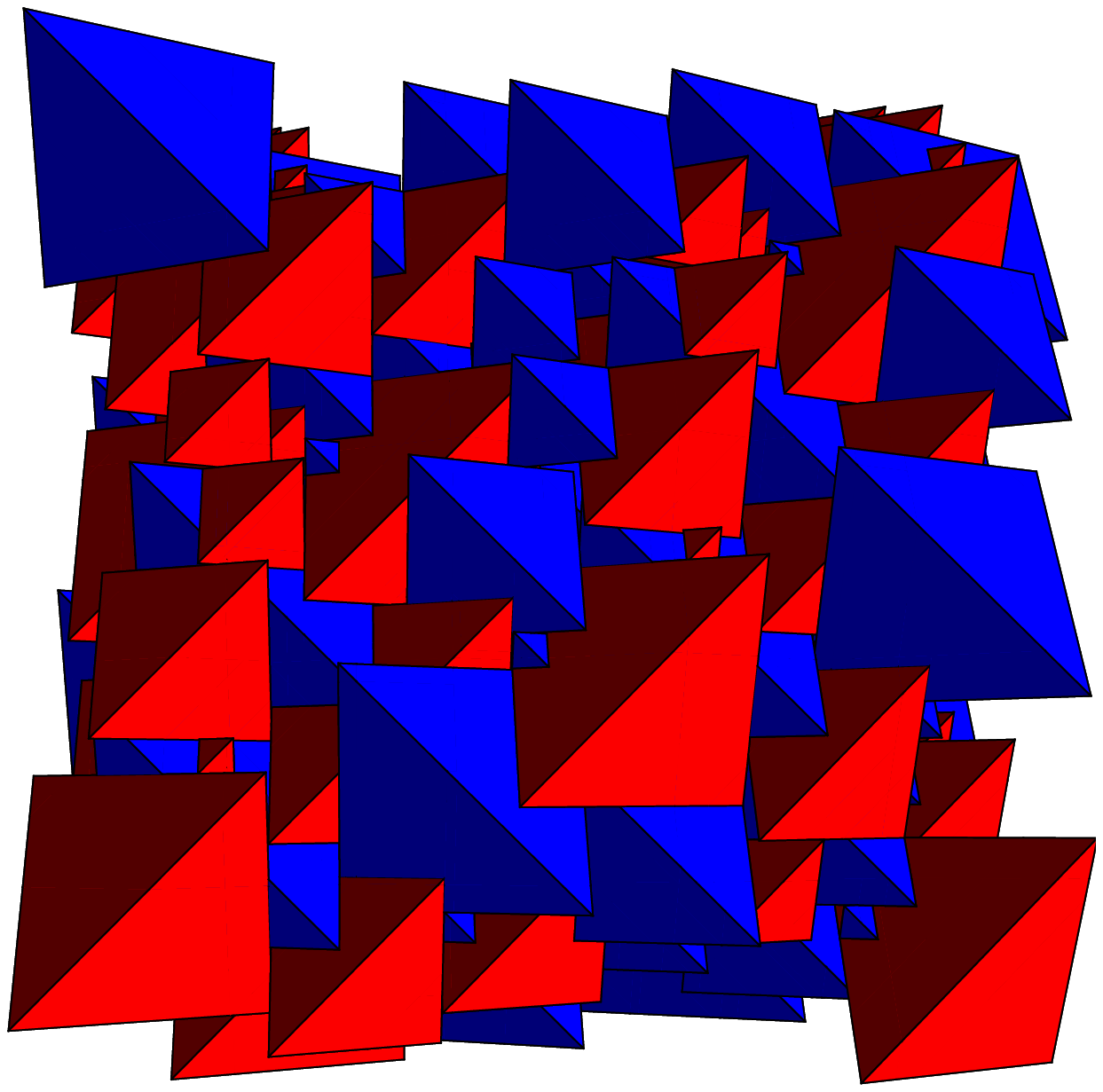}

\caption{\label{fig4}The top figure shows a two-dimensional cut through the densest known
packing of tetrahedra. The plane of
the cut is parallel to the bases of the bipyramidal dimers. Triangular sections
from dimers of one orientation (red) and from dimers of inverted orientation (blue) are
visible. The bottom figure shows a three-dimensional cut through the
densest known packing of pentatopes. The cut is taken parallel to the bases of the pentatope
dimers, and tetrahedral sections from the two dimer orientations (red and blue, again) are
visible.  
}
\end{figure}

For the problem of packing regular four-dimensional simplices (pentatopes) in four-dimensional
Euclidean space, we report a new packing discovered by our search method (Figure \ref{fig4}).
This packing, with
density $\phi=128/219\approx0.5845$, is, to our knowledge, denser than any previously reported
packing of regular pentatopes. Like the densest known tetrahedron packing,
this packing also takes the form of a double lattice of dimers (a dimer here is the union of two
cell-sharing pentatopes). This
structure, composed of a repeating unit of two oppositely oriented dimers,
repeatedly came up as the densest
in \textit{de novo} PDC searches with $p=4$ and $p=8$ pentatopes in the unit cell, whereas searches
with intermediate values of $p$ yielded sparser packings. We subsequently refined
the packing with a restricted search where the dimer was taken as the basic particle.

Note that the density reported is slightly lower than that
of the densest known packing of four-dimensional spheres 
($\phi=\pi^2/16\approx0.6169$). It remains to be determined whether this is the case because 
the optimal packing density of pentatopes is smaller than that of spheres or because
the dimer double lattice is suboptimal.
The vertex coordinates of the four primitive pentatopes and the
generating matrix of the lattice are given in Table \ref{tab4}.

\begin{table}
\centering
\begin{tabular}{|l|l|}
\hline 
primitive pentatopes & $K_1=\operatorname{conv}\{\mathbf{r}_1,\mathbf{r}_2,\mathbf{r}_3,\mathbf{r}_4,\mathbf{r}_5\}$ \\
& $K_2=\operatorname{conv}\{\mathbf{r}_2,\mathbf{r}_3,\mathbf{r}_4,\mathbf{r}_5\,\mathbf{r}_6\}$ \\
& $K_3=\mathbf{t}-K_1$\\
& $K_4=\mathbf{t}-K_2$\\
\hfill where & $\mathbf{r}_1=\sqrt{5}(1,1,1,1)$\\
& $\mathbf{r}_2=(3,-1,-1,-1)$\\
& $\mathbf{r}_3=(-1,3,-1,-1)$\\
& $\mathbf{r}_4=(-1,-1,3,-1)$\\
& $\mathbf{r}_5=(-1,-1,-1,3)$\\
& $\mathbf{r}_6=-\sqrt{5}(1,1,1,1)$\\
& $\mathbf{t}=\frac{1}{4}(-7,1,3,3)-\frac{\sqrt{5}}{4}(1,1,1,1)$\\
\hline
lattice & $\Lambda=\mathbb{Z}^4\mathrm{M}_0$\\
\hfill where& $\mathrm{M}_0 = \frac{1}{4}\left(\begin{array}{cccc}
-6 & 10 & -6 & 2 \\
-8 & -4 & 4 & 8 \\
-7 & 5 & 9 & -7 \\
1 & -7 & 9 & -3
\end{array}\right)$\\ 
&$~~+\frac{\sqrt{5}}{4}\left(\begin{array}{cccc}
2 & 2 & 2 & 2 \\
2 & 2 & 2 & 2 \\
1 & 1 & 1 & 1 \\
3 & 3 & 3 & 3
\end{array}\right)$\\
\hline

\end{tabular}

\caption{\label{tab4}Coordinates of the densest pentatope packing discovered by the PDC search ($\phi=4\operatorname{vol}(K_1)/\det(\mathrm{M}_0)=128/219\approx0.5845$).}

\end{table}

\section{\label{sec5}Conclusion}

In this article we report on the development of PDC, a novel, constraint-based method for
discovering dense periodic packings through \textit{de novo} numerical searches.
We lay out the principles of the method and demonstrate its application for selected problems. 
In addition to the dense packing of regular tetrahedra reported in Ref. \cite{tetra},
we also discover a new dense packing of regular pentatopes using the PDC method.
We also use the method to numerically recover the lattice sphere packings of highest
known density and highest known kissing number in a range of dimensions, providing
empirical evidence of their optimality.

In developing the PDC scheme, we adapt the $D-C$ framework to periodic systems. PDC
retains the mindset of the traditional $D-C$ approach of Ref. \cite{D-C}, but generalizes its
formalism in a few ways. We introduce an expanded configuration space
parameterized by linear combinations of the original parameters, such that these new parameters
over-determine the configuration. Therefore, by contrast with the traditional
construction, where new parameters are, specifically,
redundant copies of original parameters and concurrence
is described by the equality of all copies of a given original parameter, here we allow concurrence
to be described by a general linear relation. With this generalization, we can treat
the periodic images of a particle as ``replicas'' of the particle, even as they are 
related by a lattice vector instead of being identical. Thus, the 
variables describing the periodic repetition of
the configuration, namely the lattice vectors, are not imposed as constants or 
adjusted in dedicated steps. Instead, due to the projection formulation of the dynamics,
the unit cell variables that minimize the change to the configuration are determined at
each iteration. These variables  are treated on the same footing as
particle positions and orientations and are optimized as aggressively. 

Additionally, we develop a displacement-minimizing overlap resolution algorithm 
for the convex hulls of two sets of points in $\mathbb{R}^d$.
We use this algorithm to implement the projection 
to the exclusion constraint in the case of polytopal particles. 

Unlike Monte Carlo simulations, which explore the physical optimization landscape using
stochastic moves, a PDC search uses a deterministic map in an expanded, non-physical
configuration space. As such, it is useful when interest lies more in discovering optimal configurations
and less in discovering the physical pathways to such configurations. However, introducing
non-physical dynamics has been observed to be important in overcoming dynamical stagnation
\cite{Glotzer}. The projection-based
dynamics make PDC particularly well-suited in problems with hard constraints, such as hard particle
packing, or with step potentials, which prohibit the use of gradient information.

While no direct comparison has been made between the performance of PDC and Monte Carlo
searches in the case of 
periodic packing problems, difference map and $D-C$ methods in the case of other problems
have been shown to perform better than or on a par with specialized and general-purpose methods
\cite{DM-PNAS, D-C, Gravel, Protein}. The generality of the PDC scheme and its demonstrated ability
to discover dense packings in a variety of settings indicate its utility as a general method
for conducting \textit{de novo} numerical searches and as a 
possibly attractive alternative to conventional methods \footnote{An implementation of our algorithm
is available upon request from the corresponding author.}. 

Y. K. acknowledges N. Duane Loh for valuable discussions. This work was supported by grant NSF-DMR-0426568.

\bibliography{revpdc}

\end{document}